\documentclass[10pt,reqno]{amsart}

\usepackage{algorithmic}
\usepackage{algorithm}
\usepackage{amsbsy}
\usepackage{amsfonts}
\usepackage{amsmath}
\usepackage{amssymb}
\usepackage{bbm}
\usepackage{constants}
\usepackage{subcaption}
\usepackage{cite}
\usepackage{xcolor}
\usepackage{enumerate}
\usepackage{framed}     
\usepackage{graphicx}
\usepackage{multirow} 
\usepackage{tabu,longtable,supertabular,float,lipsum}
\usepackage{tabularx,url,amsmath}
\usepackage{booktabs,threeparttable}
\usepackage{pict2e}

\usepackage{pdflscape}

\usepackage[pagewise]{lineno}  
\nolinenumbers

\usepackage{hyperref}
\hypersetup{colorlinks=true,
            linkcolor=black,
            citecolor=blue,
            pdftitle={Paths in Markoff Graphs Over Finite Fields},
            pdfauthor={Silverman}}

\newcommand{\drawnfrom}{\xleftarrow{\$}}

\newtheorem{theorem}{Theorem}[section]

\newtheorem{proposition}[theorem]{Proposition}

\newtheorem{heuristic}[theorem]{Heuristic Assumption}

\theoremstyle{definition}
\newtheorem{definition}[theorem]{Definition}

\newtheorem{remark}[theorem]{Remark}

\theoremstyle{remark}

  {\end{list}}

  {\end{list}}

\renewcommand{\a}{\alpha}

\renewcommand{\b}{\beta}

\newcommand{\g}{\gamma}

\newcommand{\f}{\varphi}

\renewcommand{\l}{\lambda}

\renewcommand{\r}{\rho}

\newcommand{\s}{\sigma}

\newcommand{\z}{\zeta}

\newcommand{\Ccal}{{\mathcal C}}

\newcommand{\Mcal}{{\mathcal M}}

\newcommand{\Rcal}{{\mathcal R}}

\newcommand{\Wcal}{{\mathcal W}}

\renewcommand{\AA}{\mathbb{A}}

\newcommand{\FF}{\mathbb{F}}
\newcommand{\GG}{\mathbb{G}}

\newcommand{\PP}{\mathbb{P}}

\newcommand{\RR}{\mathbb{R}}

\newcommand{\ZZ}{\mathbb{Z}}

\newcommand{\Aut}{\operatorname{Aut}}

\newcommand{\GL}{\operatorname{GL}}

\newcommand{\prim}{\textup{prim}}
\newcommand{\Prob}{\operatornamewithlimits{Prob}}

\renewcommand{\setminus}{\smallsetminus}

\newcommand{\SL}{\operatorname{SL}}

\newcommand{\PATH}{\operatorname{\textup{\textsf{PATH}}}}
\newcommand{\FACTOR}{\operatorname{\textup{\textsf{FACTOR}}}}
\newcommand{\DLP}{\operatorname{\textup{\textsf{DLP}}}}

\newcommand{\NEGLIGIBLE}{\operatorname{\textup{\textsf{NEGLIGIBLE}}}}
\newcommand{\GEN}{\operatorname{\textsf{Gen}}}

\newcommand{\MarkoffPathFinder}{\operatorname{\textsf{MarkoffPathFinder}}}

\newcommand{\MaximalEllipticQ}{\operatorname{\textsf{MaximalEllipticQ}}}
\newcommand{\MarkoffDLP}{\operatorname{\textsf{MarkoffDLP}}}
\newcommand{\ClassicalDLP}{\operatorname{\textsf{ClassicalDLP}}}
\newcommand{\FALSEX}{\textbf{false}}  
\newcommand{\TRUEX}{\textbf{true}}
\newcommand{\COMMENTX}{\textbf{comment}:\enspace}

\newcommand{\SmallMatrix}[1]{%
  \left(\begin{smallmatrix} #1 \end{smallmatrix}\right)}

\newcommand{\TABT}[1]{\begin{tabular}[t]{@{}l@{}}#1\end{tabular}}
\newcommand{\TAB}[1]{\begin{tabular}{@{}l@{}}#1\end{tabular}}

\newcommand{\elliptic}{hyperbolic}
\newcommand{\ellipticity}{hyperbolicity}
\newcommand{\Elliptic}{Hyperbolic}
\newcommand{\hyperbolic}{elliptic}
\newcommand{\hyperbolicity}{ellipticity}

\date{\today}

\title[Paths in Markoff Graphs]
      {A Heuristic Subexponential Algorithm to Find Paths in Markoff Graphs Over Finite Fields}
\date{\today}
\author[J.H. Silverman]{Joseph H. Silverman}
\email{joseph\_silverman@math.brown.edu}
\address{Mathematics Department, Box 1917
  Brown University, Providence, RI 02912 USA.
  ORCID: 0000-0003-3887-3248}

\thanks{Silverman's research supported by Simons Collaboration Grant \#712332}

\begin{document}

\begin{abstract}
Charles, Goren, and Lauter~[J.\ Cryptology 22(1),~2009] explained how
one can construct hash functions using expander graphs in which it is
hard to find paths between specified vertices.  The set of solutions
to the classical Markoff equation $X^2+Y^2+Z^2=XYZ$ in a finite
field~$\FF_q$ has a natural structure as a tri-partite graph using
three non-commuting polynomial automorphisms to connect the
points. These graphs conjecturally form an expander family, and Fuchs,
Lauter, Litman, and Tran~[Mathematical~Cryptology~1(1),~2022]
suggested using this family of Markoff graphs in the CGL construction.
In this note we show that in both a theoretical and a practical sense,
assuming two randomness hypotheses, one can compute paths in a Markoff
graph over~$\FF_q$ by factoring~$q-1$ and solving three discrete
logarithm problems in~$\FF_q^*$. In particular, the path problem can
be solved in subexponential time.
\end{abstract}

\maketitle

\setcounter{tocdepth}{3}

\tableofcontents

\section{Introduction}
\label{section:introduction}

The classical Markoff surface is the affine surface given by the
equation
\begin{equation}
  \label{eqn:marksurf}
  \Mcal : X^2 + Y^2 + Z^2 = 3XYZ.
\end{equation}
There are three double covers $\Mcal\to\AA^2$ that give rise to three
non-commuting involutions~$\s_1,\s_2,\s_3$. A famous theorem of
Markoff~\cite{Markoff} says that every positive integer solution
of~\eqref{eqn:marksurf} can be obtained from~$(1,1,1)$ by repeatedly applying
the~$\s_i$ and permuting the coordinates.
\par
In this note we consider solutions to~\eqref{eqn:marksurf} in a finite
field~$\FF_q$ of characteristic at least~$5$. There has recently been
a lot of interest in studying the orbit structure
of~$\Mcal(\FF_q)$~\cite{BGS1,BGS,Chen,MR4112680,MR4074547,MR4213160}.
Baragar~\cite{MR2686830} conjectured that the reduction
modulo~$q$ map
\[
\Mcal(\ZZ)\longrightarrow\Mcal(\FF_q)
\]
is surjective for all primes~$q$. A recent deep result of  William Chen~\cite{Chen},
building on ground-breaking work of Bourgain, Gamburd, and Sarnak~\cite{BGS1,BGS},
says that Baragar's conjecture is true for all sufficiently large primes.\footnote{In
an updated version of~\cite{BGS1} (private communication), the authors
note that the algorithm in the present article ``should allow one to
check Baragar's conjecture for much larger~$q$. Whether it is feasible
to bridge the gap and verify the conjecture for all primes is an
interesting question.''}
\par
More precisely, consider the three non-commuting automorphisms
\[
\rho_1,\rho_2,\rho_3 : \Mcal\longrightarrow\Mcal
\]
given by the formulas
\[
\rho_1 = (X,Z,3XZ-Y),\quad
\rho_2 = (3XY-Z,Y,X),\quad
\rho_3 = (Y,3YZ-X,Z),
\]
where~$\rho_1,\rho_2,\rho_3$ are obtained from~$\s_1,\s_2,\s_3$ by
composing with appropriate coordinate permutations.
(The~$\rho_i$ are called ``rotations'' in~\cite{BGS1,BGS}.)
Let
\[
\Rcal=\langle\rho_1,\rho_2,\rho_3\rangle \subset \Aut(\Mcal)
\]
be the group of automorphisms generated by the~$\rho_i$, and let
\[
\Mcal^*(\FF_q) = \Mcal(\FF_q) \setminus \bigl\{(0,0,0)\bigr\}.
\]
Then Chen's theorem says that~$\Mcal^*(\FF_q)$ consists of a single
$\Rcal$-orbit for sufficiently large prime values of~$q$.
\par
We consider the undirected function graph (sometimes called a Schrier graph)
associated to the action of~$\{\rho_1,\rho_2,\rho_3\}$
on the set~$\Mcal^*(\FF_q)$, i.e., we form an undirected
graph~$\overline{\Mcal}(\FF_q)$ whose vertices and edges are
given by
\begin{align*}
\textsf{Vertices}\Bigl(\overline\Mcal(\FF_q)\Bigr) &= \Mcal^*(\FF_q),\\
\textsf{Edges}\Bigl(\overline\Mcal(\FF_q)\Bigr) &=
\Bigl\{ [P,\rho_i(P)] : P\in\Mcal^*(\FF_q),\; i=1,2,3 \Bigr\}.
\end{align*}
It is conjectured in~\cite{BGS1,BGS} that~$\overline\Mcal(\FF_q)$ is a
family of expander graphs; see also~\cite{MR4213160}.
\par
Charles, Goren, and Lauter~\cite{CGL} have explained how one can build
cryptographic hash functions from expander graphs provided that it
is hard to find paths in the graph connecting two given vertices.
This led Fuchs, Lauter, Litman, and Tran~\cite{FLLT} to
suggest using the Markoff graph~$\overline{\Mcal}(\FF_q)$ to construct
a hash function.  They prove, using the connectivity ideas
from~\cite{BGS1,BGS}, that there is a path-finding algorithm
for~$\overline\Mcal(\FF_q)$ that runs in deterministic
time~$O(q\log\log{q})$, and they speculate that any path-finding
algorithm in~$\overline\Mcal(\FF_q)$ must take time at least~$O(q)$.
This leads them to suggest that ``these graphs may be good candidates''
for the~CGL~hash function construction.
\par
Our goal in this note is to show that under some reasonable heuristic
assumptions, it is possible to solve the path-finding problem
in~$\overline\Mcal(\FF_q)$ in subexponential time on a classical computer
and polynomial time on a quantum computer. More precisely, up to
small polynomial-time tasks, 
it suffices to factor~$q-1$ and 
solve three discrete logarithm problems in~$\FF_q^*$, as
described in the following theorem, whose
proof will be given in Section~$\ref{section:pathfindingalgorithm}$.

\begin{theorem}[Markoff Path-Finding Algorithm] 
\label{theorem:PATHq}
We set the following notation\textup:
\begin{align*}
  \PATH(q) &= \text{time to find a path between points in $\overline\Mcal(\FF_q)$.} \\
  \DLP(q) &= \text{time to solve the \textup{DLP} in $\FF_q^*$.} \\
  \FACTOR(N) &= \text{time to factor $N$.} \\
  \NEGLIGIBLE(q) &= \TABT{tasks that take negligible (polylog) time as a function
    of~$q$,\\ for example taking square roots in $\FF_q$, or iterations\\
    performed $(q-1)/\f(q-1)\le2\log\log(q)$ times.\\}
\end{align*}  
Assume that Heuristics~$\ref{heuristic:probmaxfiber}$
and~$\ref{heuristic:probcurve}$ are valid.
Then with high probability,
\[
\PATH(q) \le \FACTOR(q-1) + 3\cdot\DLP(q) + \NEGLIGIBLE(q).
\]
\end{theorem}

\begin{remark}
We note that the paths constructed by our Markoff path-finder
algorithm (Algorithm~\ref{algorithm:markoffpathfinder})
have the following two properties that are unlikely to be
present in the paths generated by the~CGL graph-theoretic hash
function~\cite{CGL}:
\begin{itemize}
\item
  They are quite long, in the sense that the number of~$\rho_i$ used to
  connect~$P$ to~$Q$ is almost certainly larger than, say,~$q^{1/2}$.
\item
  The path connecting~$P$ to~$Q$ has long stretches in which it
  repeatedly applies one of the~$\rho_i$, e.g., it is almost certainly
  true that somewhere in the path there is a~$\rho_k$ that is repeated
  at least~$q^{1/2}$ times without using the other two~$\rho_i$.
\end{itemize}
Thus one might still consider using the Markoff graph for a CGL hash
function with the proviso that long or repetitive paths are
disallowed. On the other hand, the fact that one can create paths and
collisions, even of a disallowed type, may cause some disquiet, as well
as making it more difficult to construct a security reduction proof.
\end{remark}


\begin{remark}
In this article we have restricted attention to the classical Markoff
equation~\eqref{eqn:marksurf}, but we note that the method works, mutatis mutandis,
for more general Markoff--Hurwitz type equations
\begin{multline}  
  \label{eqn:markoffhurwitz}
  a_1X^2+a_2Y^2+a_3Z^2+b_1XY+b_2XZ+b_3YZ \\
  +c_1X+c_2Y+c_2Z+dXYZ+e=0
\end{multline}
that admit three non-commuting involutions.
\end{remark}

We give an initial high-level description of the Markoff path-finder
algorithm in Section~\ref{section:informalpathfindingalgorithm} using
psudo-code (Table~\ref{table:pathfindingalgorithm}) and a picture
(Table~\ref{table:illustratepathalg}). A more detailed description
that includes the path-finder algorithm
(Algorithm~\ref{algorithm:markoffpathfinder}) and its subroutines is
given in Section~\ref{section:pathfinderandsubroutines}.  The proof
that the Markoff path-finder algorithm finds a path and has the
indicated running time is given in
Section~\ref{section:pathfindingalgorithm}. 
A key observation in constructing the path-finder algorithm, as already exploited
in~\cite{BGS1,BGS}, is to note that the action of the~$\rho_i$ on appropriate
fibers of~$\Mcal\to\AA^1$ is described via repeated application of a
linear transformation in~$\SL_2$. (In fancier terminology, the fibers
are~$\GG_m$-torsors.) This means that if we are given two points on a
fiber, then finding a power of~$\rho_i$ that links the given points can be
rephrased as a discrete logarithm problem, either in~$\FF_q^*$ or in
the subgroup of norm~$1$ elements of~$\FF_{q^2}$. 
\par
The heuristic part of our algorithm comes from assuming that if we
take a random point in~\text{$\Mcal(\FF_q)$} or on a related curve, and if we
take one of the coordinates~$t\in\FF_q$ of that point, then
\begin{equation}
  \label{eqn:ProbT23tT1primrt}
  \Prob\left(\begin{tabular}{@{}l@{}}
    $T^2-3tT+1$ has a root\\
     $\l\in\FF_q^*$ that generates $\FF_q^*$\\
  \end{tabular}\right)
  \approx \frac12\cdot\frac{\f(q-1)}{q-1}.
\end{equation}
The factor of~$\frac12$ in~\eqref{eqn:ProbT23tT1primrt} comes from the
probability that a random quadratic polynomial has its roots~$\FF_q$,
and the factor of~$\frac{\f(q-1)}{q-1}$ comes from the well-known fact
that there are~$\f(q-1)$ generators (primitive roots) in~$\FF_q^*$.
We describe the heuristic assumptions more precisely in
Section~\ref{section:heuristic}, and we do some numerical experiments
to test the assumptions in Section~\ref{section:checkheuristic}.
\par
We illustrate the Markoff path-finding algorithm in
Section~\ref{section:pathfinderexample} by executing it on a numerical
example with $q=70687$. We find paths between some randomly chosen
points in~$\Mcal(\FF_q)$, and a non-trivial loop from a point back to
itself. Finally, in Section~\ref{section:MarkoffK3} we briefly discuss
a family of~K3 surfaces that is analogous to the Markoff
surface~\eqref{eqn:marksurf} and its
generalizations~\eqref{eqn:markoffhurwitz} and explain how
path-finding on the associated graphs can be heuristically reduced to
the elliptic curve discrete logarithm problem (ECDLP).

\par\vspace{10pt}
\noindent\textbf{Acknowledgements}.\enspace
The author would like to thank Elena Fuchs and Igor Shparlinski for their helpful comments. 

\section{A High-Level Description of the Markoff Path-Finding Algorithm}
\label{section:informalpathfindingalgorithm}

For the convenience of the reader,
Table~\ref{table:pathfindingalgorithm} gives an informal description
of our heuristically subexponential algorithm for finding paths
in~$\overline\Mcal(\FF_q)$. In this algorithm, we say that an
element~$t\in\FF_q^*$ is \emph{maximally {\elliptic}} if the quadratic
polynomial $T^2-3tT+1$ has a root~$\l\in\FF_q$ that is a primitive
root, i.e., such that~$\l$ is a generator for~$\FF_q^*$;
cf.\ Definition~\ref{definition:maximallyelliptic}.  The algorithm is
also illustrated by the picture in
Table~\ref{table:illustratepathalg}.  We refer the reader to
Section~\ref{section:pathfinderandsubroutines} for a more detailed
description of the Markoff path-finding algorithm, and to
Section~\ref{section:pathfindingalgorithm} for a proof that the
Markoff path-finding algorithm operates successfully in subexponential
time.

\begin{table}[ht]
  \def\ITEM[#1]{\item[\textup{\textbf{#1}\enspace}]}
  \framebox{\vbox{ 
\begin{itemize}
\ITEM[\textbullet]
  \textbf{Input}:
  \begin{tabular}[t]{cl}
    $\FF_q$ &  a finite field of characteristic at least~$5$ \\
    $P,Q$ & points in $\Mcal^*(\FF_q)$ \\
  \end{tabular}
\ITEM[\textbullet]
    Use $\rho_1$ and $\rho_3$ to randomly move $P$
    in $\Mcal(\FF_q)$ until reaching a point $P'$
    satisfying $y(P')$ is maximally {\elliptic}.
    This gives~$i_1,\ldots,i_\a\in\{1,3\}$ such that
    \[
    P' = \rho_{i_{\a}}\circ\cdots\circ\rho_{i_{1}}(P).
    \]
\ITEM[\textbullet]
    Use $\rho_1^{-1}$ and $\rho_2^{-1}$ to randomly move $Q$
    in $\Mcal(\FF_q)$ until reaching a point $Q'$
    satisfying $z(Q')$ is maximally {\elliptic}.
    This gives~$j_1,\ldots,j_\b\in\{1,2\}$ such that
    \[
    Q = \rho_{j_{1}}\circ\cdots\circ\rho_{j_{\b}}(Q').
    \]
\ITEM[\textbullet]
  Let $F(X,Y,Z)=X^2+Y^2+Z^2-3XYZ$. Randomly select maximally
  {\elliptic}~$x_0\in\FF_q^*$ until finding a value for which the quadratic equations
  \[
  F\bigl( x_0, y(P'), Z \bigr) =  F\bigl( x_0, Y, z(Q') \bigr) = 0
  \]
  have a solution~$(y_0,z_0)\in\FF_q^2$. Set
  \[
  P''\gets \bigl( x_0, y(P'), z_0 \bigr)
  \quad\text{and}\quad
  Q''\gets \bigl( x_0, y_0, z(Q') \bigr). 
  \]
  We note that:
  \begin{itemize}
  \item[\textbullet]
    $P''$ and~$Q''$ are on the same maximally {\elliptic}~$x$-fiber,
  \item[\textbullet]
    $P'$ and~$P''$ are on the same maximally {\elliptic}~$y$-fiber,
  \item[\textbullet]
    $Q'$ and~$Q''$ are on the same maximally {\elliptic}~$z$-fiber.
  \end{itemize}
\ITEM[\textbullet]
  Find~$a,b,c$ satisfying
  \[
  P'' = \rho_2^a(P'),\quad
  Q' = \rho_3^b(Q''),\quad
  Q'' = \rho_1^c(P'').
  \]
  As explained in Proposition~\ref{proposition:DLPonxfiber}, this
  involves solving three DLPs in~$\FF_q^*$.
\ITEM[\textbullet]
  \textbf{Output}:
  The list of integers $(i_1,\ldots,i_\a),\;(j_1,\ldots,j_\b),\;(a,b,c)$
  specifies the path
  \[
  Q = 
  \rho_{j_{1}}\circ\cdots\circ\rho_{j_{\b}}\circ
  \rho_3^b\circ
  \rho_1^c\circ
  \rho_2^a\circ
  \rho_{i_{\a}}\circ\cdots\circ\rho_{i_{1}}(P).
  \]
\end{itemize}
}}
\caption{High-level description of the Markoff path-finding algorithm}
\label{table:pathfindingalgorithm}
\end{table}

\begin{table}[ht]
  \begin{picture}(300,150)(-10,-20)
    \put(0,0){\circle*5}
    \put(5,-5){\makebox(0,0)[l]{$P''$}}
    \put(200,100){\circle*5}
    \put(205,95){\makebox(0,0)[lt]{$Q''$}}      
    \put(0,-20){\line(0,1){140}}
    \put(-20,-10){\line(2,1){240}}
    \put(100,100){\line(1,0){200}}
    \put(-5,50){\makebox(0,0)[r]{$P'$}}
    \put(0,50){\circle*5}
    \put(250,105){\makebox(0,0)[lb]{$Q'$}}
    \put(250,100){\circle*5}
    \put(125,55){\makebox(0,0)[l]{$x$-fiber}}
    \put(5,115){\makebox(0,0)[l]{$y$-fiber}}
    \put(135,105){\makebox(0,0)[b]{$z$-fiber}}

    \put(40,60){\circle*5}
    \put(43,60){\makebox(0,0)[l]{$P$}}
    \Vector(36,59)(4,51)
    \put(220,50){\circle*5}
    \put(225,50){\makebox(0,0)[l]{$Q$}}
    \Vector(221,55)(248,97)
    \put(105,33){\makebox(0,0)[lt]{\small
        \framebox{
          \begin{tabular}{@{}l@{}}
            Short random walks from $P$ to $P'$\\
            and $Q$ to $Q'$, together with three\\
            DLP computations to find paths \\
            $P'\to P''$, $P''\to Q''$, and $Q''\to Q'$\\
          \end{tabular}
          }
        }}
  \end{picture}
  \caption{Illustrating the Markoff path-finding algorithm}
  \label{table:illustratepathalg}
\end{table}

\section{Rotations on a Fiber and an Associated Matrix}
\label{section:rotationsonfiber}

The map~$\rho_1$ may be written in matrix form as
\begin{equation}
\label{eqn:iterationasmatrix}
\rho_1(x,y,z) = 
\begin{pmatrix} 1 & 0 & 0 \\ 0 & 3x & -1\\ 0 & 1 & 0 \\ \end{pmatrix}
\begin{pmatrix} x\\ y\\ z\\ \end{pmatrix}.
\end{equation}
Thus computing~$\rho_1^n(x,y,z)$ amounts to taking
the~$n$th power of the matrix~$\SmallMatrix{3x&-1\\1&0\\}$.
Similar considerations apply to~$\rho_2$ and~$\rho_3$.
This prompts the following definitions.

\begin{definition}
For~$t\in\FF_q^*$, we  set the following notation:
\[
  L_t = \begin{pmatrix} 3t & -1\\ 1 & 0 \\ \end{pmatrix} \in \SL_2(\FF_q),\qquad
  \l_t^{\vphantom1},\l_t^{-1} =\text{the eigenvalues of $L_t$.} 
\]
We note that~$\l_t\in\FF_{q^2}^*$, and that~$\l_t$ is in~$\FF_q^*$
if and only if~$9t^2-4$ is a square in~$\FF_q^*$.
\end{definition}


Formula~\eqref{eqn:iterationasmatrix} tells us that if we apply
iterates of~$\rho_1$ to a point~$(x,y,z)\in\Mcal^*(\FF_q)$, then
\begin{equation}
  \label{eqn:rho1nxyzxynzn}
  \rho_1^n(x,y,z) = (x,y_n,z_n)
  \quad\text{with}\quad
  \begin{pmatrix} y_n\\ z_n\\ \end{pmatrix}
  = L_x^n\begin{pmatrix} y\\ z\\ \end{pmatrix},
\end{equation}
and similarly for~$\rho_2$ and~$\rho_3$. This often allows us to find paths
in fibers of~$\overline\Mcal(\FF_q)$ by solving a DLP in~$\FF_q^*$.

\begin{definition}
\label{definition:probelliptic}
In~\textup{\cite{BGS1,BGS,FLLT}}, the various~$L_t$ are separated into
three cases, analogous to the classification of elements
of~$\SL_2(\RR)$.
We say that~$t\in\FF_q^*$ is
\begin{center}
\begin{tabular}{r@{:\quad}l}
  \emph{{\elliptic}} & if $\l_t\in\FF_q\setminus\{\pm1\}$,\\
  \emph{parabolic} & if $\l_t=\pm1$, \\
  \emph{{\hyperbolic}} & if $\l_t\in\FF_{q^2}^*\setminus\FF_q^*$. \\
\end{tabular}
\end{center}
\end{definition}

\begin{remark}
The characteristic polynomial of~$L_t$ is $T^2-3tT+1$, whose discriminant
is~\text{$9t^2-4$}, so we see that
\begin{align*}
  \text{$L_t$ is {\elliptic}}
  &\quad\Longleftrightarrow\quad 9t^2-4 \in {\FF_q^*}^2, \\
  \text{$L_t$ is parabolic}
  &\quad\Longleftrightarrow\quad 9t^2-4=0, \\
  \text{$L_t$ is {\hyperbolic}}
  &\quad\Longleftrightarrow\quad 9t^2-4 \notin {\FF_q^*}^2.
\end{align*}
Hence if~$q$ is large, then a randomly chosen~$t\in\FF_q$ has a roughly~$50\%$
chance of being {\elliptic}, a~$50\%$ chance of being {\hyperbolic}, and negligible
chance of being parabolic.
\end{remark}

We state a counting result for the number of~$\FF_q$-rational
points on a fiber of~$\Mcal$ and sketch the elementary proof.

\begin{proposition}
\label{proposition:numberpointsfiberM}
Let~$x_0\in\FF_q^*$. Then
\[
\#\Bigl\{ (y,z)\in\FF_q^2 : (x_0,y,z)\in\Mcal(\FF_q) \Bigr\}
= \begin{cases}
  q-1 &\text{if $x_0$ is {\elliptic},} \\
  q+1 &\text{if $x_0$ is {\hyperbolic}.} \\
\end{cases}
\]
\end{proposition}
\begin{proof}
This is proven in~\cite[Lemmas~4 and~5]{BGS1}, but
for the convenience of the reader, we sketch the proof.  
If~$x_0$ is {\elliptic} or {\hyperbolic}, i.e., if~$9x_0^2\ne4$, then 
\[
C_{x_0} : U^2 - 3x_0UV + V^2 = -x_0^2W^2
\]
is a non-singular conic in~$\PP^2$. Therefore~$C_{x_0}\cong_{/\FF_q}\PP^1$,
and hence~$\#C_{x_0}(\FF_q)=\#\PP^1(\FF_q)=q+1$.
The two points at~$\infty$, i.e., the points with~$W=0$, are
defined over~$\FF_q$ if and only if~$9x_0^2-4$ is a square in~$\FF_q^*$, i.e., 
if and only if~$x_0$ is {\elliptic}. Thus
\[
\#\Bigl( C_{x_0}\setminus\{W=0\} \Bigr)(\FF_q)
=\begin{cases}
\#\PP^1(\FF_q) = q+1 &\text{if $x_0$ is {\hyperbolic},}\\
\#\PP^1(\FF_q)-2 = q-1 &\text{if $x_0$ is {\elliptic}.}\\
\end{cases}
\tag*{\qedhere}
\]
\end{proof}

\begin{definition}
\label{definition:maximallyelliptic}
We say that~$t\in\FF_q^*$ is \emph{maximally {\elliptic}} if any one of
the following equivalent conditions is true\textup:
\begin{itemize}
\item
  $L_t$ is {\elliptic} and has order $q-1$ in $\SL_2(\FF_q)$.
\item
  An eigenvalue~$\l_t$ of~$L_t$ is in~$\FF_q$ and generates
  the multiplicative group~$\FF_q^*$.
\item
  There is a generator~$\l\in\FF_q^*$ such that $t=\frac13(\l+\l^{-1})$.  
\end{itemize}
\end{definition}

\begin{proposition}
\label{proposition:probmaxell}
We have the following formulas\textup:
\begin{align}
  \label{eqn:tmaxelnum}
  \#\{ t\in\FF_q^* : \text{$t$ is maximally {\elliptic}} \} &= \frac{\f(q-1)}{2}. \\
  \label{eqn:tmaxelprob}
  \Prob_{t\in\FF_q^*}\bigl(\text{$t$ is maximally {\elliptic}} \bigr)
  &= \frac12\cdot\frac{\f(q-1)}{q-1}.
\end{align}
\end{proposition}

\begin{proof}
Let~$\GEN(q)\subset\FF_q^*$ denote the set of generators of~$\FF_q^*$,
and consider the map
\[
f : \FF_q^*\longrightarrow\FF_q,\quad \l\longrightarrow \frac13(\l+\l^{-1}).
\]
We want to count~$\#f\bigl(\GEN(q)\bigr)$.  The map~$f$ is exactly
$2$-to-$1$ onto its image, since~$f(\l)=f(\l^{-1})$, except for the
two points~$f(\pm1)=\pm\frac23$ that have only one pre-image. The set
of generators of~$\FF_q^*$ is invariant under inversion and does not
contain~$\pm1$, so
\[
\#f\bigl(\GEN(q) \bigr) = \frac12\#\GEN(q).
\]
The group~$\FF_q^*$ is cyclic of order~$q-1$,
so~$\#\GEN(q)=\phi(q-1)$, which gives~\eqref{eqn:tmaxelnum},
and~\eqref{eqn:tmaxelprob} is an immediate corollary.
\end{proof}

We now prove the key result that if~$x\in\FF_q$ is maximally
{\elliptic}, then~$\rho_1$ acts transitively on the~$x$-fiber of~$\Mcal$,
and that we can explicitly find $\rho_1$-paths in the~$x$-fiber by
solving a DLP in~$\FF_q^*$.

\begin{proposition}
\label{proposition:DLPonxfiber}
Let $x\in\FF_q^*$ be maximally {\elliptic}, and let 
\begin{equation}
  \label{eqn:PxyzPprimexyzprime}
  P = (x,y,z)\in\Mcal^*(\FF_q) \quad\text{and}\quad P'=(x,y',z')\in\Mcal^*(\FF_q)
\end{equation}
be any two points on the $x$-fiber of~$\Mcal$.  Then there exists
an~$n\ge0$ such that
\begin{equation}
  \label{eqn:Pprimeeqrho1nP}
  P' = \rho_1^n(P),
\end{equation}
and we can compute an exponent~$n$ satisfying~\eqref{eqn:Pprimeeqrho1nP}
by solving a quadratic equation in~$\FF_q$
and then solving a DLP in the group~$\FF_q^*$.
\textup(See Algorithm~$\ref{algorithm:markoffdlp}$
in Table~$\ref{table:markoffdlp}$ for an explicit algorithm.\textup)
\end{proposition}
\begin{proof}
The assumption that~$x$ is maximally {\elliptic}
means, by definition, that the eigenvalues~$\l_x^{\vphantom1},\l_x^{-1}$ of the matrix~$L_x$
are elements of order~$q-1$ in~$\FF_q^*$. Further, since we always assume that~$q>3$, we
know that~$\l_x\ne\pm1$. 
This allows us to diagonalize~$L_x$ working over~$\FF_q$. Explicitly,
\begin{equation}
  \label{eqn:ULxUinvdiag}
  U = \begin{pmatrix}1&-\l_x^{-1}\\-1&\l_x\\  \end{pmatrix} \in \GL_2(\FF_q)
  \quad\text{satisfies}\quad
  U L_x U^{-1} = \begin{pmatrix} \l_x & 0 \\ 0 & \l_x^{-1} \\ \end{pmatrix}.
\end{equation}
(Note that~$\l_x\ne\pm1$ implies that~$U$ is invertible, since $\det(U)=\l_x^{\vphantom1}-\l_x^{-1}$.)
\par
We first prove that~$\rho_1$ acts transitively on the~$x$-fiber of~$\Mcal$. To do this,
we characterize the integers~$m$ such that~$\rho_1^m$ fixes~$P=(x,y,z)$. Thus
\begin{align*}
\rho_1^m(P)=P
&\quad\Longleftrightarrow\quad
\begin{pmatrix} y \\ z \\ \end{pmatrix}
= L_x^m \begin{pmatrix} y \\ z \\ \end{pmatrix}
\quad\text{from~\eqref{eqn:rho1nxyzxynzn},} \\
&\quad\Longleftrightarrow\quad
U \begin{pmatrix} y \\ z \\ \end{pmatrix}
=  \begin{pmatrix} \l_x^m & 0 \\ 0 & \l_x^{-m} \\ \end{pmatrix}
U  \begin{pmatrix} y \\ z \\ \end{pmatrix}
\quad\text{from \eqref{eqn:ULxUinvdiag},} \\
&\quad\Longrightarrow\quad
\left(\TAB{
  $1$ is an eigenvalue of $L_x$, since we know that
  $\SmallMatrix{y\\z\\}\ne\SmallMatrix{0\\0\\}$ and\\
  that~$U$ is invertible, so $U\SmallMatrix{y\\z\\}\ne\SmallMatrix{0\\0\\}$
  is an eigenvector\\
  of~$\SmallMatrix{\l_x^m & 0 \\ 0 & \l_x^{-m} \\}$ with eigenvalue~$1$,\\}
\right)\\
&\quad\Longrightarrow\quad
\l_x^m=1
\quad\text{since the eigenvalues of $\SmallMatrix{\l_x^m & 0 \\ 0 & \l_x^{-m} \\}$
  are $\l_x^{\pm m}$,} \\
&\quad\Longrightarrow\quad
q-1\mid m\quad\text{since $\l_x$ has order $q-1$ in $\FF_q^*$.}
\end{align*}
This proves in particular that the~$q-1$ points
\begin{equation}
  \label{eqn:Pr1Pf12Pf1q2P}
  P,\,\r_1(P),\,\r_1^2(P),\,\r_1^3(P),\,\ldots,\,\r_1^{q-2}(P)
\end{equation}
are distinct. On the other hand, Propostion~\ref{proposition:numberpointsfiberM}
and the assumption that~$x$ is {\elliptic} tell us that the~$x$-fiber of~$\Mcal$
has exactly~$q-1$ points with coordinates in~$\FF_q$, so~\eqref{eqn:Pr1Pf12Pf1q2P}
is the complete list of such point. This completes the proof
that~$\rho_1$ acts transitively on the~$x$-fiber.
\par
Now let~$P$ and~$P'$ be points~\eqref{eqn:PxyzPprimexyzprime}
on the~$x$-fiber of~$\Mcal$. We have just proven that~$\rho_1$ acts transitively,
so we know that there exists an~$n\ge0$ such that~$P'=\rho_1^n(P)$,
and we want to describe how to compute such an~$n$.
The first step is to compute~$\l_x$, which requires solving a quadratic equation.
We then repeat our earlier calculation, 
\begin{align}
\label{eqn:Pprimeequrho1nP}
P'=\rho_1^n(P)
&\quad\Longleftrightarrow\quad
\begin{pmatrix} y' \\ z' \\ \end{pmatrix}
= L_x^n \begin{pmatrix} y \\ z \\ \end{pmatrix}
\quad\text{from~\eqref{eqn:rho1nxyzxynzn},} \notag \\
&\quad\Longleftrightarrow\quad
U \begin{pmatrix} y' \\ z' \\ \end{pmatrix}
=  \begin{pmatrix} \l_x^n & 0 \\ 0 & \l_x^{-n} \\ \end{pmatrix}
U  \begin{pmatrix} y \\ z \\ \end{pmatrix}
\quad\text{from \eqref{eqn:ULxUinvdiag}.}
\end{align}
We note that the vectors
\[
U \begin{pmatrix} y' \\ z' \\ \end{pmatrix} 
= \begin{pmatrix}  y'-\l_x^{-1}z' \\ -y'+\l_x z' \\ \end{pmatrix}
\quad\text{and}\quad
U \begin{pmatrix} y \\ z \\ \end{pmatrix} 
= \begin{pmatrix}  y-\l_x^{-1}z \\ -y+\l_x z \\ \end{pmatrix}
\]
are non-zero (since~$(y,z)\ne(0,0)$ and $(y',z')\ne(0,0)$ and~$U$ is invertible)
and that they involve only the known quantities~$\l_x,y,z,y',z'\in\FF_q$. Hence the
only unknown quantity in~\eqref{eqn:Pprimeequrho1nP} is~$n$.  Since
the vectors are non-zero, at least one of the coordinates
of~\eqref{eqn:Pprimeequrho1nP} gives an equation of the
form~$\a\l_x^n=\b$ with known~$\a,\b\in\FF_q^*$, so we can find~$n$
by solving a DLP in~$\FF_q^*$.
\end{proof}

\begin{remark}
We have focused on points of $\Mcal(\FF_q)$ that have a coordinate
that is maximally {\elliptic} for reasons of computational and
expositional simplicity.  But we note that we could also use points
with a maximally {\hyperbolic} coordinate, where~$t\in\FF_q^*$ is said to
be \emph{maximally {\hyperbolic}} if~$L_t$ has order~\text{$q+1$}
in~$\SL_2(\FF_q)$, or equivalently if the roots of the
polynomial~\text{$T^2-3tT+1$} generate the subgroup of~$\FF_{q^2}^*$
of index~$q-1$. Checking for maximal {\hyperbolicity} requires a
factorization of~\text{$q+1$}, which could be an advantage
if~\text{$q+1$} is easier than~\text{$q-1$} to factor. And using
maximal {\hyperbolic} points would more-or-less double the probability of
finding a point having a fiber on which the assocaited rotation acts
transitively.  On the other hand, we would then need to solve the~DLP
in the order~\text{$q+1$} subgroup of~$\FF_{q^2}^*$, which is more difficult than
working in~$\FF_q$, although
still a subexponential problem.  In any case, the algorithm
still requires solving three~DLPs.
\end{remark}

\section{\texorpdfstring{Checking If $t\in\FF_q^*$ Is Maximally {\Elliptic}}
  {Checking If Maximally {\Elliptic}}}
\label{section:orderFqstar}
We analyze the running time of
Algorithm~\ref{algorithm:maxelliptic} in Table~\ref{table:maxelliptic},
which checks whether a given~$t\in\FF_q^*$ is maximally {\elliptic}, i.e.,
whether the matrix~$L_t=\SmallMatrix{3t&-1\\1&0}$ has order~$q-1$
in~$\SL_2(\FF_q)$. Algorithm~\ref{algorithm:maxelliptic} is then
invoked in Steps~\ref{step:MPF9},~\ref{step:MPF10},
and~\ref{step:MPF11} of the Markoff path-finder algorithm
(Algorithm~\ref{algorithm:markoffpathfinder} in
Table~\ref{table:markoffpathfinder}).
\par
The first step in Algorithm~\ref{algorithm:maxelliptic} is to
find a non-zero root of
\begin{equation}
  \label{eqn:T23tT10inFqroot}
  T^2 - 3tT + 1 = 0
\end{equation}
in~$\FF_q$, or show there is no such root. This is done by checking if
the discriminant~$9t^2-4$ of~\eqref{eqn:T23tT10inFqroot} is a square,
and if it is, using a practical polynomial-time square-root algorithm.
\par
Assuming that the equation~\eqref{eqn:T23tT10inFqroot} has a
root~$\l\in\FF_q$, it remains to check whether~$\l$ generates~$\FF_q^*$,
i.e., whether~$\l$ is a primitive root. 
The most straightforward way to check this is to first factor~$q-1$,
\[
q-1 = \prod_{i=1}^r p_i^{e_i},
\]
which need only be done once, and then use the elementary 
fact:
\begin{equation}
  \label{eqn:primrtiffpowne1}
  \text{$\l$ is a primitive root}
  \quad\Longleftrightarrow\quad
  \l^{(q-1)/p_i}\ne1~\text{for all $1\le i\le r$.}
\end{equation}
Hence once~$q-1$ has been factored, we have:
\[
\left(\begin{tabular}{@{}l@{}}
  time to check if $t\in\FF_q$\\
  is maximally {\elliptic}\\
\end{tabular}\right)
=
\left(\begin{tabular}{@{}l@{}}
  time to compute\\
  a square root in $\FF_q$\\
\end{tabular}\right)
+
\left(\begin{tabular}{@{}l@{}}
  time to compute $r$\\
  exponentiations in $\FF_q$\\
\end{tabular}\right).
\]
Since taking square roots and doing exponentiations take practical
polynomial time, and since~$r<\log_2(q)$, the time to check if an
element of~$\FF_q^*$ is maximally {\elliptic} is negligible.

\begin{remark}
\label{remark:partiallyfactorq1}
We note that the factorization of~\text{$q-1$} is used to make it easy
to check if an element of~$\FF_q^*$ is a primitive root.
However, rather than completely factoring~\text{$q-1$}, we could instead use
Lenstra's elliptic curve factorization algorithm~\cite{MR916721} to find all
moderately small prime factors. This is very efficient, since the
running time for Lenstra's algorithm to factor an integer~$N$ depends
on the size of the smallest prime factor of~$N$. We can then use the partial factorization
to create a probabilistic primitive root algorithm that has a high success
rate. Thus for example, if~$q\approx2^{4000}$ and we use Lenstra's algorithm to find
all primes~$p<2^{100}$ that divide~$q-1$, we can consider the algorithm
\begin{equation}
  \label{eqn:probprimroot}
  \text{$\l$ is probably a primitive root}
  \quad\Longleftrightarrow\quad
  \l^{(q-1)/p}\ne1~\text{for all $p\mid q-1$, $p<2^{100}$.}
\end{equation}
The probability that~\eqref{eqn:probprimroot} misidentifies an element of~$\FF_q^*$
as a primitive root when~$q\approx2^{4000}$ is less than
\[
1 - \prod_{p\mid q-1,\,p>2^{100}} \left(1-\frac{1}{p}\right)
<
1 - \left(1 - \frac{1}{2^{100}}\right)^{40}
\approx \frac{1}{2^{94}},
\]
so the probability is negligible. And even if~\eqref{eqn:probprimroot}
returns a false positive, the path-finding algorithm can simply
restart.  Finally, we note the since factoring~\text{$q-1$} and
solving the DLP in~$\FF_q^*$ are of roughly the same order of
difficulty using the best known algorithms, the saving in only
partially factoring~\text{$q-1$} is minimal at best.
\end{remark}


\section{Heuristic Assumptions}
\label{section:heuristic}

In this section we describe the heuristic assumptions that we will use
in our Markoff path-finding algorithm. They say roughly that if we
choose a~$t\in\FF_q^*$ that is the coordinate of a random point
in~$\Mcal(\FF_q)$ or a random point on a certain curve, then the
probability that~$t$ is maximally {\elliptic} is roughly the same as
if~$t$ were chosen randomly in~$\FF_q^*$.  We refer the reader to
Section~\ref{section:checkheuristic} for data that supports
Heuristics~\ref{heuristic:probmaxfiber} and~\ref{heuristic:probcurve}.

\begin{heuristic}
\label{heuristic:probmaxfiber}
Let $(x_0,y_0,z_0)\in\Mcal(\FF_q)$. For~$n\ge0$,
randomly choose~$i_1,i_2,\ldots,i_n\in\{1,3\}$ and set
\[
(x_n,y_n,z_n)
\gets\rho_{i_n}\circ\rho_{i_{n-1}}\circ\cdots\circ\rho_{i_2}\circ\rho_{i_1}(x_0,y_0,z_0).
\]
Then
\begin{equation}
  \label{eqn:probmax2}
  \Prob \bigl(\text{$y_n$ is maximally {\elliptic}}  \bigr)
  \approx\frac12 \frac{\f(q-1)}{q-1}.
\end{equation}
\end{heuristic}
\begin{proof}[Justification]
As we randomly use~$\rho_1$ and~$\rho_3$ to ``rotate'' on the
$x$-fiber and the~$z$-fiber, it is reasonable to
view~$y_0,y_1,y_2,\ldots$ as being independent random elements
of~$\FF_q^*$, at least insofar as to whether they are maximally
{\elliptic}, which recall means that an associated quadratic polynomial
has a root that generates~$\FF_q^*$.    Proposition~\ref{proposition:probmaxell}
says that a random element of~$\FF_q^*$ has
probability~$\f(q-1)/2(q-1)$ of being maximally {\elliptic}, which gives
the desired justification.
\end{proof}

\begin{heuristic}
\label{heuristic:probcurve}
Let
\[
F(X,Y,Z)=X^2+Y^2+Z^2-3XYZ,
\]
let~$a,b\in\FF_q^*$, and let~$\Ccal_{a,b}\subset\AA^3$ be the affine
curve
\[
\Ccal_{a,b} : F\bigl( X, a, Z \bigr) =  F\bigl( X, Y, b \bigr) = 0.
\]
Then
\begin{equation}
  \label{eqn:probtisxandmaxell}
  \Prob_{t\in\FF_q^*}\left(
  \begin{tabular}{@{}l@{}}
    $t$ is maximally {\elliptic} and is the\\
    $x$-coordinate of a point in $\Ccal_{a,b}(\FF_q)$\\
  \end{tabular}
  \right) \approx \frac18\cdot\frac{\f(q-1)}{q}. 
\end{equation}
\end{heuristic}

\begin{proof}[Justification]
We note that~$t\in\FF_q^*$ is the $x$-coordinate of a point
in~$\Ccal_{a,b}(\FF_q)$ if and only if the quadratic equations
\[
F(t,a,Z)=0 \quad\text{and}\quad F(t,Y,b)=0
\]
have roots in~$\FF_q$, so if and only if
\[
9a^2t^2-4(t^2+a^2) \in {\FF_q^*}^2
\quad\text{and}\quad
9b^2t^2-4(t^2+b^2) \in {\FF_q^*}^2.
\]
Hence we expect that
\[
\Prob_{t\in\FF_q^*}\Bigl(
\text{$t$ is the $x$-coordinate of a point in $\Ccal_{a,b}(\FF_q)$}\Bigr)
\approx \frac14,
\]
since about half the elements of~$\FF_q^*$ are squares. (If~$a^2=b^2$,
the probability increases to~$\frac12$, which helps the attacker.)
\par
It is reasonable to view the~$x$-coordinates of the points
in~$\Ccal_{a,b}(\FF_q)$ as independent random elements of~$\FF_q$.
Proposition~\ref{proposition:probmaxell} says that a random element
of~$\FF_q^*$ has probability~$\frac{\f(q-1)}{2(q-1)}$ of being
maximally {\elliptic}, so multiplying this by the probability~$\frac14$
that~$t$ is the $x$-coordinate of a point in~$\Ccal_{a,b}(\FF_q)$
yields~\eqref{eqn:probtisxandmaxell}.
We also note that Weil's estimate says that
\[
\#\Ccal_{a,b}(\FF_q) = q + O(\sqrt{q}),
\]
so when~$q$ is large, the set~$\Ccal_{a,b}(\FF_q)$ is also large.
\end{proof}

\begin{remark}
\label{remark:heuristic3viaweil}
Igor Shparlinski has noted that one might be able to prove
Heuristic~\ref{heuristic:probcurve} using estimates for points on
curves over finite fields and the covering argument in~\cite{BGS1}.
We sketch the ideas in
Section~\ref{section:covermaxellpts}, although we do not give a
completely rigorous proof, which would require delicate arguments
regarding irreducibility of covers of the~$\Ccal_{a,b}$ curves.
\end{remark}

\section{The Markoff Path-Finder Algorithm}
\label{section:pathfindingalgorithm}
We give the proof of our main result (Theorem~$\ref{theorem:PATHq}$),
which we restate for the convenience of the reader.

\begin{theorem}[Theorem $\ref{theorem:PATHq}$] 
We set the following notation\textup:
\begin{align*}
  \PATH(q) &= \text{time to find a path between points in $\overline\Mcal(\FF_q)$.} \\
  \DLP(q) &= \text{time to solve the \textup{DLP} in $\FF_q^*$.} \\
  \FACTOR(N) &= \text{time to factor $N$.}\\
  \NEGLIGIBLE(q) &= \TABT{tasks that take negligible (polylog) time as a function
    of~$q$,\\ for example taking square roots in $\FF_q$, or iterations\\
    performed $(q-1)/\f(q-1)\le2\log\log(q)$ times.\\}
\end{align*}  
Assume that Heuristics~$\ref{heuristic:probmaxfiber}$
and~$\ref{heuristic:probcurve}$ are valid.  Then with high
probability, the Markoff path-finder Algorithm described in detail as
Algorithm~$\ref{algorithm:markoffpathfinder}$ in
Table~$\ref{table:markoffpathfinder}$ will find a path between
randomly given points in the graph~$\overline\Mcal(\FF_q)$ in time
\begin{equation}
  \label{eqn:pathfactor3dlp}
  \PATH(q) \le \FACTOR(q-1) + 3\cdot\DLP(q) + \NEGLIGIBLE(q).
\end{equation}
\end{theorem}
\begin{proof}
The Markoff path-finder algorithm
(Algorithm~\ref{algorithm:markoffpathfinder} in
Table~\ref{table:markoffpathfinder}) terminates with a list of
positive integers
\[
(i_1,\ldots,i_\a),\;(j_1,\ldots,j_\b),\;(a,b,c)
\]
satisfying
\begin{align*}
  P' &= \rho_{i_{\a}}\circ\rho_{i_{\a-1}}\circ\cdots\circ\rho_{i_{2}}\circ\rho_{i_{1}}(P)
  &&\text{Steps \ref{step:MPF1}--\ref{step:MPF2}} \\
  Q &= \rho_{j_{1}}\circ\rho_{j_{2}}\circ\cdots\circ\rho_{j_{\b-1}}\circ\rho_{j_{\b}}(Q')
  &&\text{Steps \ref{step:MPF3}--\ref{step:MPF4}} \\
  P'' &= \rho_2^a(P'),\;Q'=\rho_3^b(Q''),\;Q''=\rho_1^c(P'')
  &&\text{Steps \ref{step:MPF7}--\ref{step:MPF8}} 
\end{align*}
We use these to compute
\begin{align*}
Q
&= \rho_{j_{1}}\circ\cdots\circ\rho_{j_{\b}} (Q') \\
&= \rho_{j_{1}}\circ\cdots\circ\rho_{j_{\b}} \circ \rho_3^{b} (Q'') \\
&= \rho_{j_{1}}\circ\cdots\circ\rho_{j_{\b}} \circ \rho_3^{b} \circ \rho_1^c (P'') \\
&= \rho_{j_{1}}\circ\cdots\circ\rho_{j_{\b}} \circ \rho_3^{b} \circ \rho_1^c \circ \rho_2^a (P') \\
&= \rho_{j_{1}}\circ\cdots\circ\rho_{j_{\b}} \circ \rho_3^{b} \circ \rho_1^c \circ \rho_2^a
\circ\rho_{i_{\a}}\circ\cdots\circ\rho_{i_{1}}(P).
\end{align*}
Hence Algorithm~\ref{algorithm:markoffpathfinder} gives the a path
in~$\overline\Mcal(\FF_q)$ from~$P$ to~$Q$.
\par
We next consider the running time of each step of the algorithm.  In
Step~\ref{step:MPF12} we factor the integer~$q-1$.  Once this is done,
the time to check whether an element~$t\in\FF_q$ is maximally {\elliptic}
is negligible; see Section~\ref{section:orderFqstar}.
\par
In Steps~\ref{step:MPF1}--\ref{step:MPF2} and
Steps~\ref{step:MPF3}--\ref{step:MPF4}, we randomly move a  point
around~$\Mcal(\FF_q)$ and check whether one of its coordinates is
maximally {\elliptic}.  Heuristic~\ref{heuristic:probmaxfiber} says that
each of these loops needs to look at an average of~$2(q-1)/\f(q-1)$ points
before terminating, and as already noted, checking maximal {\ellipticity} 
takes negligible time once we have factored~$q-1$.  Similarly,
Heuristic~\ref{heuristic:probcurve} says that the loop in
Steps~\ref{step:MPF11} is executed an average of~$8(q-1)/\f(q-1)$
times, with the maximal {\ellipticity} and the square root
computations taking negligible time. Hence
Steps~\ref{step:MPF1}--\ref{step:MPF6} take average time~$12(q-1)/\f(q-1)$
multiplied by some small power of~$\log(q)$.
There are classical estimates~\cite[Sections~18.4 and~22.9]{HW6}
\[
\frac{N}{\f(N)}\le \Cl{HW2}\cdot\log\log N\quad\text{for all $N\ge5$,}
\]
and indeed one can take~$\Cr{HW2}=2$, which allows us to conclude that
Steps~\ref{step:MPF1}--\ref{step:MPF6} take a neglible amount of time. 
\par
Steps~\ref{step:MPF7}--\ref{step:MPF8} use the $\MarkoffDLP$ algorithm
three times, and the $\MarkoffDLP$ algorithm
(Algorithm~\ref{algorithm:markoffdlp} in Table~\ref{table:markoffdlp})
requires taking a square root in~$\FF_q^*$ (neglible time) and
computing a discrete logarithm in~$\FF_q^*$. Hence the time to execute
Steps~\ref{step:MPF7}--\ref{step:MPF8} is essentially the time
it takes to compute three DLPs in~$\FF_q^*$.
\par
Adding these time estimates yields~\eqref{eqn:pathfactor3dlp},
which completes the proof that the Markoff path-finding algorithm
terminates in the specified time.
\end{proof}

\section{The Markoff Path Finder Algorithm in Action: An Example}
\label{section:pathfinderexample}

We illustrate the Markoff path-finder algorithm
(Algorithm~\ref{algorithm:markoffpathfinder} in Table~\ref{table:markoffpathfinder})
by computing a numerical
example. We take
\begin{align*}
  q &= 70687,\qquad q-1 = 2\cdot 3^{3}\cdot 7\cdot 11\cdot 17, \\
  P &= (45506, 13064, 18) \in\Mcal(\FF_q),\\
  Q &= (11229, 5772, 56858) \in\Mcal(\FF_q).
\end{align*}
We use a simplified version of the algorithm in which $i_\a=1$ for all~$\a$
and~$j_\b=1$ for all~$\b$, since in practice this almost always works.
Thus Steps~\ref{step:MPF1}--\ref{step:MPF2} say to apply~$\rho_1$
to~$P$ until the~$y$-coordinate is maximally {\elliptic}. We do a similar
computation in Steps~\ref{step:MPF3}--\ref{step:MPF4}, except now we
apply iterates of~$\rho_1^{-1}$ to~$Q$ and stop when we reach an
iterate whose~$z$-coordinate is maximally {\elliptic}.
Table~\ref{table:examplePiterates}
show our computations. It lists the iterates and indicates whether the
appropriate coordinate is {\elliptic} or {\hyperbolic}; and if the coordinate is
{\elliptic}, it lists~$o(\lambda)$, the order of an associated eigenvalue
in~$\FF_q^*$.  We find that
\[
y\bigl(\rho_1^2(P)\bigr)
\quad\text{and}\quad
z\bigl(\rho_1^{-15}(Q)\bigr)
\quad\text{are maximally {\elliptic},}
\]
so the output from Steps~\ref{step:MPF1}--\ref{step:MPF4} are
\begin{align*}
\a&=2, &i_1=i_2&=1,&
P' &= \rho_1^2(P) = (45506, 40902, 10340),\\
\b&=15, &j_1=\cdots=j_{15}&=1, &
Q' &= \rho_1^{-15}(Q) = (11229, 2424, 19535).
\end{align*}

\begin{table}[t]
  \[
  \begin{array}{|r@{}l|r@{}l|r@{}l|l|} \hline
  i &= 0 & P &= (45506, 13064, 18) & y &= 13064 & \text{{\elliptic}, $o(\lambda)=1683$}  \\ \hline
  i &= 1 & \rho_1(P) &= (45506, 18, 40902) & y &= 18
  & \text{{\elliptic}, $o(\lambda)=4158$}  \\ \hline
  i &= 2 & \rho_1^{2}(P) &= (45506, 40902, 10340) & y &= 40902
  & \text{{\elliptic}, $o(\lambda)=70686$}  \\ \hline
  \multicolumn{4}{c}{\hfill} \\ \hline
  j &= 0 & Q &= (11229, 5772, 56858) & z &= 56858 & \text{{\hyperbolic}}  \\ \hline
  j &= 1 & \rho_1^{-1}(Q) &= (11229, 65943, 5772) & z &= 5772
  & \text{{\elliptic}, $o(\lambda)=35343$}  \\ \hline
  j &= 2 & \rho_1^{-2}(Q) &= (11229, 6407, 65943) & z &= 65943
  & \text{{\elliptic}, $o(\lambda)=10098$}  \\ \hline
  j &= 3 & \rho_1^{-3}(Q) &= (11229, 29942, 6407) & z &= 6407 & \text{{\hyperbolic}}  \\ \hline
  j &= 4 & \rho_1^{-4}(Q) &= (11229, 16944, 29942) & z &= 29942 & \text{{\hyperbolic}}  \\ \hline
  j &= 5 & \rho_1^{-5}(Q) &= (11229, 35748, 16944) & z &= 16944 & \text{{\hyperbolic}}  \\ \hline
  j &= 6 & \rho_1^{-6}(Q) &= (11229, 2200, 35748) & z &= 35748 & \text{{\hyperbolic}}  \\ \hline
  j &= 7 & \rho_1^{-7}(Q) &= (11229, 66363, 2200) & z &= 2200 & \text{{\hyperbolic}}  \\ \hline
  j &= 8 & \rho_1^{-8}(Q) &= (11229, 21119, 66363) & z &= 66363 & \text{{\hyperbolic}}  \\ \hline
  j &= 9 & \rho_1^{-9}(Q) &= (11229, 46109, 21119) & z &= 21119 & \text{{\hyperbolic}}  \\ \hline
  j &= 10 & \rho_1^{-10}(Q) &= (11229, 47313, 46109) & z &= 46109
  & \text{{\elliptic}, $o(\lambda)=594$}  \\ \hline
  j &= 11 & \rho_1^{-11}(Q) &= (11229, 7133, 47313) & z &= 47313 & \text{{\hyperbolic}}  \\ \hline
  j &= 12 & \rho_1^{-12}(Q) &= (11229, 47632, 7133) & z &= 7133
  & \text{{\elliptic}, $o(\lambda)=5049$}  \\ \hline
  j &= 13 & \rho_1^{-13}(Q) &= (11229, 47838, 47632) & z &= 47632 & \text{{\hyperbolic}}  \\ \hline
  j &= 14 & \rho_1^{-14}(Q) &= (11229, 19535, 47838) & z &= 47838
  & \text{{\elliptic}, $o(\lambda)=7854$}  \\ \hline
  j &= 15 & \rho_1^{-15}(Q) &= (11229, 2424, 19535) & z &= 19535
  & \text{{\elliptic}, $o(\lambda)=70686$}  \\ \hline
  \end{array}
  \]
  \caption{$\rho_1$ iterates of~$P$  until reaching a maximally {\elliptic} $y$-fiber
  and $\rho_1^{-1}$ iterates of~$Q$  until reaching a maximally {\elliptic} $z$-fiber}
  \label{table:examplePiterates}
\end{table}

In Steps~\ref{step:MPF5}--\ref{step:MPF6} we randomly
choose~$x\in\FF_q^*$ and check whether~$x$ is maximally {\elliptic} and
whether the quadratic equations
\[
F(x,40902,Z) = 0 \quad\text{and}\quad F(x,Y,19535) = 0
\]
have solutions~$y,z\in\FF_q$. It took~$5$ tries, as listed in
Table~\ref{table:step3}. So the output from
Steps~\ref{step:MPF5}--\ref{step:MPF6} consists of the two points
\[
P'' = (29896,40902,935)
\quad\text{and}\quad
Q'' = (29896,595,19535).
\]

\begin{table}
\[
\begin{array}{|c|c|c|c|} \hline
  x & F(x,40902,Z) & F(x,Y,19535) \\ \hline\hline
  29628 & \text{irreducible} & \text{irreducible} \\\hline
  19562 & \text{irreducible} & (Y-42621)(Y-57310) \\\hline
  43036 & \text{irreducible} & \text{irreducible} \\\hline
  6057  & (Z-27506)(Z-70305) & \text{irreducible} \\\hline
  29896 & (Z-935)(Z-45089) & (Y-595)(Y-6503) \\\hline
\end{array}
\]
\caption{Finding a point on $F(x,40902,Z)=F(x,Y,19535)=0$}
\label{table:step3}
\end{table}

In Steps~\ref{step:MPF7}--\ref{step:MPF8} we find a path on the
$y$-fiber from~$P'$ to~$P''$, a path on the~$z$-fiber from~$Q''$
to~$P'$, and a path on the~$x$ fiber from~$Q''$ to~$P''$. This is
done using the Markoff DLP Algorithm
(Algorithm~\ref{algorithm:markoffdlp} in
Table~\ref{table:markoffdlp}), which uses
Proposition~\ref{proposition:DLPonxfiber} to convert the path problem
in a maximal {\elliptic} fiber into a discrete logarithm
problem in~$\FF_q^*$.  Implementing this algorithm, we find that
\[
P'' = \rho_2^{26986}(P'),\quad
Q' = \rho_3^{65193}(Q''),\quad
Q'' = \rho_1^{30287}(P'').
\]
\par
Finally, the algorithm outputs
\[
(1,1),\,(1,1,\ldots,1),\,(26986, 30287, 65193),
\]
where the second item is a~$15$-tuple.
We check that this gives a path from~$P$ to~$Q$ by computing
\begin{align*}
P &= (45506, 13064, 18) \\
\rho_1^{2}(P) &= (45506, 40902, 10340) \\
\rho_2^{26986}\circ\rho_1^{2}(P) &= (29896, 40902, 935) \\
\rho_1^{30287}\circ\rho_2^{26986}\circ\rho_1^{2}(P) &= (29896, 595, 19535) \\
\rho_3^{65193}\circ\rho_1^{30287}\circ\rho_2^{26986}\circ\rho_1^{2}(P) &= (11229, 2424, 19535) \\
\rho_1^{15}\circ\rho_3^{65193}\circ\rho_1^{30287}\circ\rho_2^{26986}\circ\rho_1^{2}(P)
&= (11229, 5772, 56858) = Q.
\end{align*}
\par
If we run the algorithm a second time, the randomness in
Steps~\ref{step:MPF5}--\ref{step:MPF6} means that we are likely to
obtain a different path. (And if we hadn't simplified the choices of the~$i_\a$ and~$j_\b$,
that randomness would also lead to different paths.)
For example, using the same~$(q,P,Q)$ as
input and running the algorithm again, we obtained the output
\[
(a,c,b) = ( 26703, 52102, 29583), 
\]
which gives the path
\[
Q = \rho_1^{15}\circ\rho_3^{29583}\circ\rho_1^{52102}\circ\rho_2^{26703}\circ\rho_1^{2}(P).
\]
\par
We also note that we can combine a path from~$P$ to~$Q$ with a path
from~$Q$ to~$P$ to 
find a non-trivial loop that starts and
returns to~$P$, since it is unlikely that the two paths will be exact
inverses of one another. Indeed, running the algorithm to find a path from~$Q$ to~$P$,
we found
\[
(1,1,1),\,\underbrace{(1,1,\ldots,1)}_{\text{$11$-tuple}},\,(389,14491,39906),
\]
which gives the path
\[
P = \rho_1^{11}\circ\rho_3^{39906}\circ\rho_1^{14491}\circ\rho_2^{389}\circ\rho_1^{3}(Q).
\]
Combining this with the first path from~$P$ to~$Q$ that we found earlier gives the loop
\[
P =
\rho_1^{11}\circ\rho_3^{39906}\circ\rho_1^{14491}\circ\rho_2^{389}
\circ\rho_1^{18}\circ\rho_3^{65193}\circ\rho_1^{30287}\circ\rho_2^{26986}\circ\rho_1^{2}(P)
\]
where we have combined the middle~$\rho_1^3\circ\rho_1^{15}$ into a
single~$\rho_1^{18}$.

\section{Markoff-Type K3 Surfaces and the ECDLP}
\label{section:MarkoffK3}

In this section we briefly discuss~K3 surfaces that are analogous to
the Markoff surface. These surfaces, which were dubbed tri-involutive~K3 (TIK3) surfaces
in~\cite{arxiv2201.12588}, are surfaces
\[
\Wcal \subset \PP^1\times\PP^1\times\PP^1
\]
given by the vanishing of a~$(2,2,2)$ form. With appropriate
non-degeneracy conditions, the three double
covers~\text{$\Wcal\to\PP^1\times\PP^1$} give three non-commuting
involutions~$\s_1,\s_2,\s_3\in\Aut(\Wcal)$. If the~$(2,2,2)$-form is
symmetric, then~$\Wcal$ also admits coordinate permutation
automorphisms, in which case we can define the analogues of the
rotations~$\rho_1,\rho_2,\rho_3\in\Aut(\Wcal)$. Fuchs, Litman, Tran, and the
present author studied the orbit structure of~$\Wcal(\FF_q)$ for
various groups of autormorphism.  In view of~\cite{FLLT}, one might
consider using the graph structure on~$\Wcal(\FF_q)$ induced
byt~$\{\s_1,\s_2,\s_3\}$ or~$\{\rho_1,\rho_2,\rho_3\}$ to implement
the CGL~\cite{CGL} hash function algorithm. However, the three
fibrations $\Wcal\to\PP^1$ have genus~$1$ fibers, the Jacobians of
these fibrations are elliptic surfaces of rank at least~$1$, and the
action of the automorphisms on fibers can be described in terms of
translation by a section to the elliptic surface. See for
example~\cite{MR3519436}, where this geometry is explained and
explicit formulas are provided.
\par
Thus the Markoff path-finder algorithm, with suitable tweaks,
yields a~K3~path-finder algorithm whose running time is determined
primarily by how long it takes to solve three instances of the
elliptic curve discrete logarithm problem. Thus on a classical
computer, the algorithm currently takes exponential time to find paths
in~$\Wcal(\FF_q)$, but that is reduced to polynomial time on a quantum
computer. However, since the algorithm can look at many elliptic
curves lying in the fibration~\text{$\Wcal(\FF_q)\to\PP^1(\FF_q)$}, it
may well be possible to find one whose order is fairly smooth, in
which case the ECDLP becomes easier to solve.  We have not pursued
this further, but it might be interesting to see if under reasonable
heuristic assumptions, one can solve the path-finding problem
in~$\Wcal(\FF_q)$ in subexponential time on a classical computer.

\appendix

\section{Sketch of Proof to Rigorously Justify Heuristic \ref{heuristic:probcurve}}
\label{section:covermaxellpts}
In this section we sketch how one might prove Heuristic
\ref{heuristic:probcurve} using Weil's estimate and an
inclusion-exclusion calculation modeled after an argument described
in~\cite{BGS1}.
Let~$\Ccal/\FF_q$ be a smooth projective curve of genus~$g$, and
let~$x:\Ccal\to\PP^1$ be a coordinate function on~$\Ccal$. We want to estimate
\[
\#\bigl\{ P\in\Ccal(\FF_q) : \text{$x(P)$ is maximally {\elliptic}} \bigr\}.
\]
We consider the following double cover of~$\Ccal$,
\[
\Ccal' = \bigl\{ (P,\l) \in C\times\PP^1 : \l^2-3x(P)\l+1=0 \bigr\},
\]
so we want to estimate the size of
\[
\Ccal'(\FF_q)^\prim :=
\bigl\{ (P,\l)\in\Ccal'(\FF_q) : \text{$\l$ is a primitive root in $\FF_q^*$} \bigr\}.
\]
\par
For each~$n\mid{q-1}$ we define 
\[
\Ccal'[n] = \bigl\{ (P,\mu)\in C\times\PP^1 : (P,\mu^n) \in \Ccal' \bigr\}.
\]
We will assume that the curves~$\Ccal'[n]$ are irreducible.  (We note
that experimentally this appears to be true for the
curves~$\Ccal_{a,b}$ appearing in Heuristic \ref{heuristic:probcurve},
but a rigorous proof would require non-trivial ideas.)
\par
The map
\[
F_n : \Ccal'[n]\longrightarrow\Ccal',\quad
F_n(P,\mu) = (P,\mu^n)
\]
is an $n$-fold cover of~$\Ccal'$, and up to a negligible number of points,
the induced map
\begin{equation}
  \label{eqn:Fnnto1onFqpts}
  F_n : \Ccal'[n](\FF_q)\longrightarrow\Ccal'(\FF_q)
  \quad\text{is also $n$-to-$1$.}
\end{equation}
\par
We will also assume that all curves are smooth and projective, which
simplifies our calculations, since we will ignore 
sets of points of negligible size.  The Riemann--Hurwitz formula tells us that the genera
of the~$\Ccal'[n]$ satisfy
\[
g\bigl(\Ccal'[n]\bigr)
= 1 + n \bigl(g(\Ccal')-1\bigr) + \frac12 \sum_{\g\in\Ccal'} ( e_\g-1 )
= O\bigl( n \bigr),
\]
where the big-$O$ constant depends only on the genus of~$\Ccal$ and
the degree of the map~$x:\Ccal\to\PP^1$.  Then Weil's estimate for the
number of points on curves over finite fields yields
\begin{equation}
  \label{eqn:numptCprime}
  \#\Ccal'[n](\FF_q) = q + O\bigl(g\bigl(\Ccal'[n]\bigr)\sqrt{q}\bigr)
  = q + O(n\sqrt{q}).
\end{equation}
\par
We note that~$(P,\l)\in\Ccal'(\FF_q)$ satisfies:
\begin{equation}
  \label{eqn:Plorderdivq1n}
  (P,\l)\in F_n\bigl( \Ccal'[n](\FF_q) \bigr)
  \quad\Longleftrightarrow\quad
  \text{the order of $\l$ in $\FF_q^*$ divides $\dfrac{q-1}{n}$.}
\end{equation}
Ignoring a small number of points that are singular or ``at
infinity,'', we use this to calculate
\begin{align}
  \label{eqn:CpFqprimnum}
  \#\Ccal'(\FF_q)^\prim
  &\approx
  \sum_{n\mid q-1} \mu\left(\frac{q-1}{n}\right) \#F_n\bigl(\Ccal'[n](\FF_q)\bigr)
  \quad\text{inclusion/exclusion and \eqref{eqn:Plorderdivq1n},}
  \notag\\
  &\approx
  \sum_{n\mid q-1} \mu\left(\frac{q-1}{n}\right)\cdot \frac{1}{n}\#\Ccal'[n](\FF_q)
  \quad\text{from \eqref{eqn:Fnnto1onFqpts},} \notag\\
  &= 
  \sum_{n\mid q-1} \mu\left(\frac{q-1}{n}\right)\cdot \left(\frac{q}{n} + O(\sqrt{q})\right)
  \quad\text{from \eqref{eqn:numptCprime},}  \notag\\
  &= q\cdot\frac{\f(q-1)}{q-1} + O\bigl(d(q-1)\sqrt{q}\bigr),
\end{align}
where~$\f$ is Euler's phi function and~$d(N)$ is the number of
divisors of~$N$. We note that~$d(N)$ is small compared to~$N$, so that
factor may be ignored.
The map~$\Ccal'\to\Ccal$ is $2$-to-$1$, so we obtain
\begin{align*}
  \#\biggl( t\in\FF_q^* \biggm|
  &
  \begin{tabular}{l@{}}
    $t$ is maximally {\elliptic} and\\
    equals $x(P)$ for some $P\in\Ccal(\FF_q)$\\
  \end{tabular}
  \biggr) \\
  &\gtrapprox
  \frac{1}{\deg x} \#\bigl\{P\in\Ccal(\FF_q) : \text{$x(P)$ is maximally {\elliptic}} \bigr\} \\
  &\gtrapprox
  \frac{1}{2\deg x} \#\Ccal'(\FF_q)^\prim \\
  &= \f(q-1) + (\text{lower order term})
  \quad\text{from \eqref{eqn:CpFqprimnum}.}
\end{align*}
This completes our sketch of how one might go about rigorously proving
Heuristic~\ref{heuristic:probcurve}.

\section{Computations to Check Heuristics
  \ref{heuristic:probmaxfiber} and \ref{heuristic:probcurve}}
\label{section:checkheuristic}

\begin{remark}[\textbf{Testing Heuristic~\ref{heuristic:probmaxfiber}}]
We choose a random point~$P$ in~$\Mcal(\FF_q)$ and randomly
apply~$\rho_1$ or~$\rho_3$ until the~$y$-coordinate of the resulting
point is maximally {\elliptic}.  For each prime in
Table~\ref{table:heur1}, we compute the average value of~$n$
for~$10^5$ randomly chosen points. We compare this with the
theoretical value~$2(q-1)/\f(q-1)$, which is the theoretical expected
number of trials to find a maximally {\elliptic} element in~$\FF_q^*$.
\end{remark}

\begin{remark}
We note that the experimental values in Table~\ref{table:heur1} are
somewhat larger than expected, especially when~$q-1$ is quite
smooth. We are not sure what is causing the discrepency. It is
probably not due to~$\Mcal(\FF_q)$ having a fewer than expected number
of points with a maximally {\elliptic} coordinate, since we checked 
experimentally in Table~\ref{table:heur1a} that the proportion of
points~$P\in\Mcal(\FF_p)$ such that~$x_P$ is maximally {\elliptic} is
almost identical to the proportion of elements of~$\FF_q$ that are
maximally {\elliptic}.  In any case, the experimental numbers are small
enough that even for~$q$ of cryptographic size, the number of
iterations of Steps~\ref{step:MPF1}--\ref{step:MPF2} and
Steps~\ref{step:MPF3}--\ref{step:MPF4} in the Markoff path-finder
algorithm (Algorithm~\ref{algorithm:markoffpathfinder} in
Table~\ref{table:markoffpathfinder}) will be practical.
\end{remark}

\begin{table}
\[
\begin{array}{|c|c|c|c|c|} \hline
  q & q-1 & 2\cdot\dfrac{q-1}{\f(q-1)}
  & \begin{tabular}{@{}l@{}} \text{Heuristic \ref{heuristic:probmaxfiber}}\\
      \text{Experiment}\\ \end{tabular} \\
  \hline\hline
  17389 & 2^{2}\cdot 3^{3}\cdot 7\cdot 23 & 7.318 & 8.919 \\ \hline
  48611 & 2\cdot 5\cdot 4861 & 5.001 & 5.200 \\ \hline
  55163 & 2\cdot 27581 & 4.000 & 3.634 \\ \hline
  70687 & 2\cdot 3^{3}\cdot 7\cdot 11\cdot 17 & 8.181 & 10.459 \\ \hline
  104729 & 2^{3}\cdot 13\cdot 19\cdot 53 & 4.662 & 4.613 \\ \hline
  200560490131 & 2\cdot3\cdot5\cdots\cdot29\cdot31 & 13.085 & 19.361 \\ \hline
\end{array}
\]
\caption{Experiments to test Heuristic \ref{heuristic:probmaxfiber} (100000 samples)}
\label{table:heur1}
\end{table}

\begin{table}
\[
\begin{array}{|c|c|c|c|}\hline
  q & \dfrac{\#\Mcal(\FF_q)}{q^2}
  &  \dfrac{\#\Mcal(\FF_q)_x^{\textup{gen}}}{q^2}   
  & \dfrac{\f(q-1)}{2(q-1)} \\ \hline\hline
  647 & 1.00000 & 0.22699 & 0.22291 \\ \hline
  757 & 1.00794 & 0.14298 & 0.14286 \\ \hline
  863 & 1.00000 & 0.24930 & 0.24942 \\ \hline
  983 & 1.00000 & 0.25034 & 0.24949 \\ \hline
  1091 & 1.00000 & 0.19853 & 0.19817 \\ \hline
  1213 & 1.00495 & 0.16514 & 0.16502 \\ \hline
  1307 & 1.00000 & 0.25003 & 0.24962 \\ \hline
\end{array}
\]
\caption{Experiments to test the size of
  $\Mcal(\FF_q)_x^{\textup{gen}}$, which is the number of $P\in\Mcal(\FF_q)$
  such that $x_P$ is maximally {\elliptic} in~$\FF_q^*$. (100000 samples)}
\label{table:heur1a}
\end{table}

\begin{remark}[\textbf{Testing Heuristic~\ref{heuristic:probcurve}}]
For each prime in Table~\ref{table:heur2}, we
chose~$10^5$ random values~$t,a,b\in\FF_q$ and
checked to see if~$t$ was both maximally {\elliptic} and
had the property that there are~$y,z\in\FF_q$
satisfying~$F(t,a,z)=F(t,y,b)=0$. We computed the proportion of
such~$t$ values and compared it with the
theoretical value~$\f(q-1)/8(q-1)$. The theoretical and
experimental values are in quite good agreement.
\end{remark}

\begin{table}
\[
\begin{array}{|c|c|c|c|c|} \hline
  q & q-1 & \dfrac{1}{8}\cdot\dfrac{\f(q-1)}{q-1}
  & \begin{tabular}{@{}l@{}} \text{Heuristic \ref{heuristic:probcurve}}\\
      \text{Experiment}\\ \end{tabular} \\
  \hline\hline
  17389 & 2^{2}\cdot 3^{3}\cdot 7\cdot 23 & 0.0342 & 0.0343 \\ \hline
  48611 & 2\cdot 5\cdot 4861 & 0.0500 & 0.0498 \\ \hline
  55163 & 2\cdot 27581 & 0.0625 & 0.0629 \\ \hline
  70687 & 2\cdot 3^{3}\cdot 7\cdot 11\cdot 17 & 0.0306 & 0.0314 \\ \hline
  104729 & 2^{3}\cdot 13\cdot 19\cdot 53 & 0.0536 & 0.0550 \\ \hline
  200560490131 & 2\cdot3\cdot5\cdots\cdot29\cdot31 & 0.0191 & 0.0192 \\ \hline
\end{array}
\]
\caption{Experiments to test Heuristic \ref{heuristic:probcurve} (100000 Samples)}
\label{table:heur2}
\end{table}

\section{The Markoff Path-Finder Algorithm and Subroutines}
\label{section:pathfinderandsubroutines}

The following algorithms are described in
Tables~\ref{table:markoffpathfinder}--\ref{table:maxelliptic} in this section.
\begin{description}
\item[Algorithm~\ref{algorithm:markoffpathfinder} - $\normalfont\MarkoffPathFinder$:]
  Returns a path in $\overline\Mcal(\FF_q)$ from $P$ to $Q$.
\item[Algorithm~\ref{algorithm:markoffdlp} - $\normalfont\MarkoffDLP$:]
  Returns an integer $n\ge0$ so that $P=\rho_k^n(G)$ in
  $\Mcal(\FF_q)$.
\item[Algorithm~\ref{algorithm:maxelliptic} -
  $\normalfont\MaximalEllipticQ$:] Returns \TRUEX\ if~$t$ is maximal
  {\elliptic} in~$\FF_q^*$, i.e., if the
  matrix~$\SmallMatrix{3t&-1\\1&0\\}$ has order~$q-1$
  in~$\SL_2(\FF_q)$; otherwise returns \FALSEX.  It assumes that a
  factorization of~$q-1$ is known; but see
  Remark~\ref{remark:partiallyfactorq1}.
\end{description}

\begin{table}[ht]
  \newcommand{\MM}[1]{}  
  \noindent
  \begin{algorithm}[H]
    \captionof{algorithm}{$\MarkoffPathFinder$}
    \label{algorithm:markoffpathfinder}
    \begin{algorithmic}[1]
      \REQUIRE $q$,\,$P$,\,$Q$ with $P,Q\in\Mcal(\FF_q)$
      \STATE \COMMENTX
      Use a factorization algorithm to factor $q-1$ and store it so that it is
      accessible by subroutines.
      \STATE $\text{PrimeFactorList}\gets\{\text{primes that divide $q-1$}\}$
      \label{step:MPF12} \MM{12}
      \STATE \COMMENTX
      Randomly move $P$ using $\rho_1$ and $\rho_3$ until the
      $y$-coordinate is maximally {\elliptic}
      \label{step:MPF1} \MM{1}
      \STATE $P'\gets P,\; \a\gets0$
      \WHILE{$\MaximalEllipticQ\bigl(q,y(P')\bigr)=\FALSE$
        \label{step:MPF9} \MM{9}
      }
      \STATE $\a\gets\a+1$
      \STATE $i_\a \drawnfrom \{1,3\}$
      \STATE $P'\gets\rho_{i_\a}(P')$
      \ENDWHILE
      \STATE \COMMENTX
      \label{step:MPF2} \MM{2}
      $P'=\rho_{i_{\a}}\circ\rho_{i_{\a-1}}\circ\cdots\circ\rho_{i_{2}}\circ\rho_{i_{1}}(P)$
      \STATE \COMMENTX
      \label{step:MPF3} \MM{3}
      Randomly move $Q$ using $\rho_1^{-1}$ and $\rho_2^{-1}$ until the
      $z$-coordinate is maximally {\elliptic}
      \STATE $Q'\gets Q,\; \b\gets0$
      \WHILE{$\MaximalEllipticQ\bigl(q,z(Q')\bigr)=\FALSE$
        \label{step:MPF10} \MM{10}
      }      
      \STATE $\b\gets\b+1$
      \STATE $j_\b \drawnfrom \{1,2\}$
      \STATE $Q'\gets\rho^{-1}_{j_\b}(Q')$
      \ENDWHILE
      \STATE \COMMENTX
      $Q=\rho_{j_{1}}\circ\rho_{j_{2}}\circ\cdots\circ\rho_{j_{\b-1}}\circ\rho_{j_{\b}}(Q')$
      \label{step:MPF4} \MM{4}
      \STATE \COMMENTX
      \label{step:MPF5} \MM{5}
      Find random points with the same maximally {\elliptic} $x$-coordinate that can be used
      to connect $P'$ to $Q'$
      \REPEAT
      \STATE $x\drawnfrom\FF_q^*$
      \UNTIL{$\MaximalEllipticQ(q,x)=\TRUE$
        \AND $F\bigl(x,y(P'),Z\bigr)$ has a root $z\in\FF_q$
        \AND $F\bigl(x,Y,z(Q')\bigr)$ has a root $y\in\FF_q$
        \label{step:MPF11} \MM{11}
      }
      \STATE $P''\gets \bigl( x, y(P'), z \bigr)$
      \STATE $Q''\gets \bigl( x, y, z(Q') \bigr)$
      \STATE\COMMENTX
      \begin{tabular}[t]{@{\textbullet\enspace}l}
        $P''$ and~$Q''$ are on the same maximally {\elliptic}~$x$-fiber.\\
        $P'$ and~$P''$ are on the same maximally {\elliptic}~$y$-fiber.\\
        $Q'$ and~$Q''$ are on the same maximally {\elliptic}~$z$-fiber.\\
      \end{tabular}
      \label{step:MPF6} \MM{6}
      \STATE \COMMENTX
      \begin{tabular}[t]{@{\textbullet\enspace}l}
      Find fibral paths $P''\to P$ and $P'\to Q''$ and $Q''\to P''$.\\
      Proposition~\ref{proposition:DLPonxfiber} ensures that such paths exist.\\
      \end{tabular}
      \label{step:MPF7} \MM{7}
      \STATE $a\gets\MarkoffDLP(q,P',P'',2)$
      \STATE $b\gets\MarkoffDLP(q,Q'',P',3)$
      \STATE $c\gets\MarkoffDLP(q,P'',Q'',1)$
      \STATE\COMMENTX
      $P''=\rho_2^a(P'),\;Q'=\rho_3^b(Q''),\;Q''=\rho_1^c(P'')$
      \label{step:MPF8} \MM{8}
      \ENSURE $(i_1,\ldots,i_\a),\;(j_1,\ldots,j_\b),\;(a,b,c)$  
    \end{algorithmic}
  \end{algorithm}
  \caption{Returns a path in $\overline\Mcal(\FF_q)$ from $P$ to $Q$}
  \label{table:markoffpathfinder}
\end{table}

\begin{table}[ht]
  \noindent
  \begin{algorithm}[H]
    \captionof{algorithm}{$\MarkoffDLP$}
    \label{algorithm:markoffdlp}
    \begin{algorithmic}[1]
      \REQUIRE $q$,\,$P$,\,$Q$,\, $k$ with $P,Q\in\Mcal(\FF_q)$ and $k\in\{1,2,3\}$
      and the $k$th coordinate of~$P$ maximally {\elliptic}
      \STATE\COMMENTX if $k=2$ ($y$-fiber) or $k=3$ ($z$-fiber), swap coordinates
      to use the $x$-fiber%
      \IF{$k=2$}
      \STATE $P\gets \bigl(y_P,z_P,x_P\bigr)$
      \STATE $Q\gets \bigl(y_Q,z_Q,x_Q\bigr)$
      \ELSIF{$k=3$}
      \STATE $P\gets \bigl(z_P,x_P,y_P\bigr)$
      \STATE $Q\gets \bigl(z_Q,x_Q,y_Q\bigr)$
      \ENDIF
      \STATE\COMMENTX
      Now $x_P$ is maximally {\elliptic}.
      \STATE $\l\gets \bigl(3x_P+\sqrt{9x_P^2-4}\bigr)/2$ in $\FF_q$
      \STATE\COMMENTX
      The maximal {\ellipticity} of $x_P$ says that $\l$ generates $\FF_q^*$.
      \STATE $b\gets \bigl( y_Q-z_Q/\l \bigr) \big/ \bigl( y_P - z_P/\l\bigr)$
      \STATE Use a DLP algorithm to find $n$ so that $\l^n=b$ in $\FF_q^*$.\\
      \ENSURE $n$
    \end{algorithmic}
  \end{algorithm}
  \caption{Returns an integer $n\ge0$ so that $P=\rho_k^n(G)$ in
    $\Mcal(\FF_q)$.  See Proposition~\ref{proposition:DLPonxfiber} for an
    explanation of why this algorithm works.}
  \label{table:markoffdlp}
\end{table}
\begin{table}[ht]
  \noindent
  \begin{minipage}[H][][t]{0.80\textwidth}
  \begin{algorithm}[H]
    \captionof{algorithm}{$\MaximalEllipticQ$} 
    \label{algorithm:maxelliptic}
    \begin{algorithmic}[1]
      \REQUIRE $q$,\,$t$
      \STATE $\textit{result}\gets\FALSE$
      \STATE\COMMENTX
      Check if $T^2-3tT+1$ has two distinct roots in $\FF_q$
      \IF{$(9t^2-4)^{(q-1)/2}=1$ in $\FF_q$}
      \STATE $\lambda\gets(3t+\sqrt{9t^2-4})/2$ in $\FF_q$
      \STATE $\textit{result}\gets\TRUE$
      \FOR{$p\in\text{PrimeFactorList}$}
      \IF{$\lambda^{(q-1)/p}=1$ in $\FF_q^*$}
      \STATE $\textit{result}\gets\FALSE$
      \ENDIF
      \ENDFOR
      \ENDIF
      \ENSURE \textit{result}      
    \end{algorithmic}
  \end{algorithm}
  \end{minipage}
  \caption{Check whether~$t$ is maximal {\elliptic}, or equivalently, whether
    the matrix $L_t\gets\SmallMatrix{3t&-1\\1&0\\}$ has order~$q-1$ in~$\SL_2(\FF_q)$, or
    equivalently, whether~$L_t$ has an eigenvalue in $\FF_q^*$ that is a primitive
    root.}
  \label{table:maxelliptic}
\end{table}

\clearpage 

\bibliographystyle{abbrv}
\bibliography{PATH}
\vspace*{20pt}

\end{document}